\documentstyle[11pt]{article}

\advance \topmargin by -\headheight
\advance \topmargin by -\headsep

\def\cros{\raise1.9pt\hbox{$\scriptscriptstyle
	   >$}\!\raise1.5pt\hbox{$\scriptstyle\triangleleft$}}

\def\o{\otimes}

\newcommand{\setc}[1]{\setcounter{equation}{#1}}
\newcounter{eqnr}

\newenvironment{eqnarrayabc}{\stepcounter{equation}
  \setcounter{eqnr}{\value{equation}}\setc{0}
  
  \begin{eqnarray}}{\end{eqnarray}\setc{\value{eqnr}}}

\renewcommand{\sec}{\setc{0}\section} 
\newcommand{\bealph}{\begin{eqnarrayabc}}
\newcommand{\eealph}{\end{eqnarrayabc}}
 
\newtheorem{prop}{Proposition}[section]
\newtheorem{lemma}[prop]{Lemma}
\newtheorem{theorem}[prop]{Theorem}
\newtheorem{defi}[prop]{Definition}

\newcommand{\be}{\begin{equation}}
\newcommand{\ee}{\end{equation}}
\newcommand{\bea}{\begin{eqnarray}}
\newcommand{\eea}{\end{eqnarray}}

\newcommand{\nn}{ \nonumber \\ }

\newcommand{\qed}{\vrule height 5pt width 5pt depth 0pt}

\newcommand{\id}{{\rm id}}
\newcommand{\1}{1\! \! 1}

\newcommand{\cp}{\raise1.9pt\hbox{$\scriptscriptstyle>$}\!\raise1.5pt
		 \hbox{$\scriptstyle\triangleleft$}}
\newcommand{\inv}{^{-1}}
\newcommand{\p}{^{\prime}}

\newcommand{\pair}[2]{\mbox{$\langle #1\vert #2 \rangle$}}

\def\cat{\mbox{$^C{\cal M}(H)_A$}}

\begin{document}

\vspace{1truecm} 
\title{\sc Doi-Hopf Modules over Weak Hopf Algebras}
 
\author{ Gabriella B\"ohm $^1$\\
\\
Research Institute for Particle and Nuclear Physics,
Budapest\\
H-1525 Budapest 114, P.O.Box 49, Hungary}
\date{}
 
\maketitle
 
\footnotetext[1]{e-mail: BGABR@rmki.kfki.hu\\
Supported by the Hungarian Scientific Research Fund, OTKA -- T 016
233}
\begin{abstract}
The theory of Doi-Hopf modules \cite{Doi,Kop} is generalized to 
Weak Hopf Algebras \cite{BSz,N,BNSz}.
\end{abstract}


\sec{Introduction}

The category $\cat$ of Doi-Hopf Modules over the bialgebra $H$ was 
introduced in \cite{Doi} and independently in \cite{Kop}. It is the 
category of the modules over the algebra $A$ which are also comodules 
over the coalgebra $C$ and
satisfy certain compatibility condition involving $H$. The study 
of $\cat$ turned out to be very useful: It was shown in 
\cite{Doi, CMZ1} that many categories investigated independently before
-- such as the module and comodule categories over bialgebras, the 
Hopf modules category \cite{Sw}, and the Yetter-Drinfeld category 
\cite{Y,RT} -- are special cases of $\cat$. Using this observation
many results known for  module categories over bialgebras or Hopf 
algebras were generalized to this more general setting \cite{CMZ2,CIMZ}.

In this paper we generalize the definition of Doi-Hopf modules to the
case when $H$ is a Weak Bialgebra (WBA). Our definitions are supported by 
the fact that many results of \cite{Kop, CMZ2, CIMZ} remain valid in this 
case.

Weak Bialgebras (Weak Hopf Algebras -- WHA's --) are generalizations of 
bialgebras (Hopf algebras) see \cite{BSz,BNSz} and \cite{N} (latter one
using somewhat different terminology). In contrast to another
direction of generalization, the quasi-Hopf algebras and weak quasi-Hopf
algebras, WBA's are coassociative. Though their counit is not an algebra
map, their structure is designed such a way that their (left or right) (co-)
module category carries a monoidal structure \cite{N, BNSz2} (and some 
more in the WHA case \cite{BNSz2}). 

WHA's have relevance for example in describing depth 2 reducible 
inclusions \cite{NSzW}. 

As the bialgebra (Hopf algebra) also the WBA (WHA) is a self-dual 
structure: The dual space of a finite dimensional WBA (WHA) carries 
naturally a WBA (WHA) structure \cite{BSz,BNSz}.

The paper is organized as follows: we define and examine the structures 
such as the weak Doi-Hopf datum  (generalizing the Doi-Hopf datum of 
\cite{Doi}) the weak Doi-Hopf module (generalizing the Doi-Hopf module of 
\cite{Doi})
the weak smash product (generalizing the analogous notion of \cite{Kop})
and the weak Doi-Hopf integral (generalizing definitions of \cite{CIMZ, 
Doint}). We illustrate these notions on the same four examples generalizing
some classical examples of \cite{ Doi, CMZ1}.

\sec{The Weak Doi-Hopf Datum}

In this Section $H$ is a Weak Bialgebra (WBA) in the sense of \cite{BNSz}
over the field $k$. Its unit element is denoted by $\1$, the product of
the elements $g,h\in H$ by $gh$, the coproduct of $h\in H$ by 
$\Delta(h)=h_{(1)}\o h_{(2)}$ and the counit is denoted by $\varepsilon$.

\begin{defi} Let $H$ be a WBA over the field $k$. The $k$-algebra $A$ is 
a {\em left $H$-comodule algebra} if there exists a left weak coaction $\rho$
of $H$ on $A$ which is also an algebra map. I.e. a map $\rho:A\to H\o A$
such that
\bealph     (\id_H\o\rho)\circ \rho &=& (\Delta\o \id_A)\circ \rho \label{a1}\\
         (\1\o a)\rho(1_A)&=&(\Pi^R\o\id_A)\circ \rho({\mathit a}) \label{a2}\\
         \rho(ab)&=&\rho(a)\rho(b) \label{a3}
\eealph
for all $a,b\in A$. We use the standard notation $\rho(a)=a_{<-1>}\o a_{<0>}$
and $(\Delta\o\id_A)\circ\rho({\mathit a})={\mathit a}_{<-2>}\o 
{\mathit a}_{<-1>}\o {\mathit a}_{<0>}= 
(\id_H\o \rho)\circ \rho ({\mathit a})$. 

The left weak coaction $\rho$ is {\em non-degenerate} if 
$(\varepsilon\o \id_A)\circ \rho= \id _A$ or, equivalently, 
$(\varepsilon \o \id_A)\circ \rho (1_A)=1_A$. 
For non-degenerate left weak coactions $\rho$ 
(\ref{a2}) has an equivalent 
form (compare with 
\cite{NSzW}) $(\Delta \o \id_A)\circ \rho(1_A)=(\1\o \rho(1_A))(\Delta(\1)\o 
1_A)$.

Similarly, $A$ is a {\em right $H$-comodule algebra} if there exists a  
right weak coaction $\rho$ of $H$ on $A$ which is also an algebra map. I.e.
a map $\rho:A\to A\o H$ such that 
\bealph (\rho\o \id_H)\circ \rho &=&(\id_A\o \Delta)\circ \rho \label{a1p}\\
        \rho(1_A)(a\o \1) &=& (\id_A \o \Pi^L)\circ\rho({\mathit a}) \label{a2p}\\
          \rho(ab) &=& \rho (a) \rho(b)  \label{a3p}
\eealph
for all $a,b\in A$. We also denote $\rho(a)=a_{<0>}\o a_{<1>}$.

The right weak coaction $\rho$ is {\em non-degenerate} 
if $(\id_A \o \varepsilon)\circ \rho= \id_A$ or, equivalently, if 
$(\id_A\o \varepsilon)\circ \rho (1_A)=1_A$.
For non-degenerate right weak coactions $\rho$
(\ref{a2p}) has an equivalent form $(\id_A\o\Delta)\circ \rho(1_A)=
(1_A\o \Delta(\1))(\rho(1_A)\o \1)$.

The dual notion to comodule algebra is the module coalgebra defined as 
follows: The $k$-coalgebra $C$ is a {\em right $H$-module coalgebra} if
there exists a right weak action of $H$ on $C$ which is also a coalgebra map. 
I.e.  a map $\cdot: C\times H\to C$ such that
\bealph (c\cdot g)\cdot h &=& c\cdot(gh) \label{c1}\\
         c\cdot \Pi^L(h) &=& \varepsilon_C(c_{(1)}\cdot h) c_{(2)} \label{c2}\\
        \Delta_C(c\cdot h)&=&\Delta_C(c)\cdot\Delta(h)  \label{c3}
\eealph 
for all $c\in C,\ g,h\in H$.  

The right weak
action $\cdot$ is {\em non-degenerate} if $c\cdot \1=c \ \forall c\in C$
or, equivalently, if $\varepsilon_C(c\cdot \1)=\varepsilon_C(c)\ \forall c\in 
C$. For non-degenerate right weak actions $\cdot$ (\ref{c2}) has the  
equivalent reformulation 
as $\varepsilon_C(c\cdot h)=\varepsilon (c\cdot \Pi^L(h))$.

Similarly, $C$ is a {\em left $H$-module coalgebra} if there exists a 
left weak action of $H$ on $C$ which is also a coalgebra map. I.e. a map
$\cdot :H\times C\to C$ such that
\bealph g\cdot (h\cdot c) &=& (gh)\cdot c \label{c1p} \\
         \Pi^R(h)\cdot c &=& c_{(1)}\varepsilon_C(h\cdot c_{(2)})\label{c2p}\\
        \Delta_C(h\cdot c)&=& \Delta(h)\cdot \Delta_C(c) \label{c3p}
\eealph
for all $c\in C\ g,h\in H$.  

The left weak action $\cdot$
is {\em non-degenerate} if $\1\cdot c=c \ \forall c\in C$ or, equivalently,
if $\varepsilon_C(\1\cdot c)=\varepsilon_C(c) \ \forall c\in C$.
For non-degenerate left weak actions $\cdot$ (\ref{c2p}) has the 
equivalent reformulation 
$\varepsilon_C (h\cdot c) = \varepsilon (\Pi^R(h)\cdot c)$.
\end{defi}

Notice, that in contrast to the case when $H$ is an ordinary bialgebra 
the unit preserving property of $\rho$ and 
the counit preserving property of $\cdot$ are  not required and the form 
of condition (b) in each group is somewhat different from the usual one.

\begin{defi} A {\em right Weak Doi-Hopf datum} is a triple $(H,A,C)$, where
$H$ is a WBA over $k$, $A$ a left $H$-comodule algebra and $C$ a right 
$H$-module coalgebra. 

A {\em left Weak Doi-Hopf datum} is a triple $(H,A,C)$ where $H$ is a WBA
over $k$, $A$ a right $H$-comodule algebra and $C$ a left $H$-module 
coalgebra.

A (left or right) Weak Doi-Hopf datum is {\em non-degenerate} if  both 
the weak coaction of
$H$ on $A$ and the weak action of $H$ on $C$ are non-degenerate. 

\end{defi}

{\noindent \bf Examples:} 

{\bf 1} Let $H$ be a WBA over $k$, $A\colon =H$ as an algebra with 
the coaction $\rho\colon =\Delta$, $C\colon = H^L$ with the coalgebra 
structure 
\bea  \Delta_{H^L}(a^L)\colon &=& \1_{(2)} a^L\o S(\1_{(1)})\equiv 
                                  \1_{(2)}\o  a^L S(\1_{(1)})\nn
      \varepsilon_{H^L}(a^L)\colon &=&\varepsilon(a^L)
\nonumber\eea
and the action $a^L\cdot h\colon = \1_{(2)}\varepsilon (a^Lh\1_{(1)})$
for all $a^L\in H^L, h\in H$. Then $(H,A=H,C=H^L)$ is a non-degenerate 
right Weak Doi-Hopf datum.

{\bf 2} Let $H$ be a WBA over $k$, $A\colon = H^L$ as the subalgebra of $H$
with the coaction $\rho\colon = \Delta\vert_{H^L}$, $C\colon = H$ as a 
coalgebra with the action 
$c\cdot h\colon = ch$ for all $c,h\in H$. Then $(H, A=H^L, C=H)$ is a 
non-degenerate right Weak Doi-Hopf datum.

{\bf 3} Let $H$ be a WBA over $k$, $A\colon = H$ as an algebra with 
the coaction $\rho\colon = \Delta$, $C\colon = H$ as a coalgebra with 
the action $c\cdot h=ch$ for all $c,h\in H$. Then $(H,A=H,C=H)$ is a 
non-degenerate right  Weak Doi-Hopf datum.

{\bf 4} Let $K$ be a WHA over $k$, $H\colon =  K^{op}\o K$ as a bialgebra.
($K^{op}$ is the bialgebra with the same coalgebra structure as $K$ and the 
opposite algebra structure.) $A\colon =K$ as an algebra with the coaction 
$\rho(a)\colon = (S\inv(a_{(3)})\o a_{(1)})\o a_{(2)}$ for all $a\in K$,
$C\colon = K$ as a coalgebra with the action $c\cdot(a\o b)\colon = acb$
for all $c\in K, (a\o b)\in H$. Then $(H=K^{op}\o K, A=K, C=K)$ is a 
non-degenerate right weak Doi-Hopf datum.
\smallskip

Let us call a (left or right) weak Doi-Hopf datum {\em finite 
dimensional} if all $H, A$ and $C$ are finite dimensional as $k$-spaces. 
There is a well defined notion of duality for finite dimensional weak Doi-
Hopf data sending a left weak Doi-Hopf datum to a right one and vice versa:

Introduce the following notations: For any finite dimensional $k$-
space $M$ let ${\hat M}$ denote the dual $k$-space. If $A$ is a finite 
dimensional algebra then by ${\hat A}$ we mean the dual space equipped with 
the dual coalgebra structure. Similarly, for a finite dimensional coalgebra 
$C$ denote the dual algebra by $\hat C$ and finally for a finite dimensional 
bialgebra $H$ denote the dual bialgebra by $\hat H$.

\begin{prop} \label{dual}
For a (non-degenerate) right weak Doi-Hopf datum $(H,A,C)$ the triple 
$({\hat H},
{\hat C}, {\hat A})$ is a (non-degenerate) left weak Doi-Hopf datum -- 
called the dual of 
$(H,A, C)$ --  with
\bea {\hat \rho}({\hat c})\colon &=& b_i\triangleright {\hat c}\o \beta^i \nn
       \phi\cdot {\hat a}\colon &=& (\phi\o {\hat a})\circ \rho 
\label{dur}
\eea
where ${\hat c}\in {\hat C}$, $\{b_i\}$ is any basis in $H$ and $\{\beta^i\}$ 
is the dual basis in ${\hat H}$, $\phi\in {\hat H}$, ${\hat a}\in {\hat A}$
and $(h\triangleright {\hat c})(d)={\hat c}(d\cdot h)$ for ${\hat c}\in 
{\hat C}, d\in C, h\in H$.

Similarly, for a (non-degenerate) left weak Doi-Hopf datum $(H,A,C)$ the 
triple $({\hat H},
{\hat C}, {\hat A})$ is a (non-degenerate) right weak Doi-Hopf datum -- 
called the dual of
$(H,A,C)$ -- with
\bea {\hat \rho}({\hat c})\colon &=&\beta^i\o  {\hat c}\triangleleft b_i\nn
     {\hat a}\cdot \phi &=&({\hat a}\o \phi)\circ \rho 
\label{dul}
\eea 
with the  obvious notation. The above duality transformation is involutive.    
\end{prop} 

{\em Proof:} The transformations (\ref{dur}) and (\ref{dul}) are
obviously inverses of one other. One easily checks that (\ref{a1p}) for 
$({\hat H},{\hat C},{\hat A})$ is equivalent to (\ref{c1}) on $(H,A,C)$, 
(\ref{a2p}) to (\ref{c2}), (\ref{a3p}) to (\ref{c3}), (\ref{c1p}) to 
(\ref{a1}), (\ref{c2p}) to (\ref{a2}) and (\ref{c3p}) to (\ref{a3}).

The non-degeneracy of the weak coaction ${\hat{\rho}}$ of $\hat{H}$ on 
${\hat C}$ 
is equivalent to the non-degeneracy of the weak action of $H$
on $C$ while the non-degeneracy of the action of $\hat{H}$ on $\hat{A}$
is equivalent to the non-degeneracy of the weak coaction $\rho$ of $H$ 
on $A$ both in the left and right cases.
\hfill\qed

\sec{The Weak Doi-Hopf Module}

\begin{defi} The $k$-space $M$ is a {\em right weak Doi-Hopf module} over the 
right weak Doi-Hopf datum $(H,A,C)$ if it is a non-degenerate right $A$-
module and a non-degenerate left $C$-comodule i.e. there exists an action 
$\cdot :M\times A\to M$ for which $m\cdot 1_A=m \ \forall m\in M$ and a 
coaction $\rho_M:M\to C\o M$ for which $(\varepsilon_C\o \id_M)\circ 
\rho_M=\rho_M$ such that  the compatibility condition 
\be \rho_M(m\cdot a)=m_{<-1>}\cdot a_{<-1>}\o m_{<0>}\cdot a_{<0>} \ee
holds for $\rho_M(m)\equiv m_{<-1>}\o m_{<0>}$.

Similarly, $M$ is a {\em left weak Doi-Hopf module} over the left Doi-Hopf
datum $(H,A,C)$ if it is a non-degenerate left $A$-module (with $A$-action 
$\cdot$) and a non-degenerate right $C$-comodule (with $C$-coaction $\rho_M$)
such that
\be \rho_M(a\cdot m)=a_{<0>}\cdot m_{<0>}\o a_{<1>}\cdot m_{<1>}.\ee
\end{defi}

The category $\cat$ has as objects the finite dimensional right weak Doi-Hopf 
modules $M$ over 
the right weak Doi-Hopf datum $(H,A,C)$ and arrows $T:M\to M\p$ which 
intertwine  both the $A$-actions and the $C$-coactions:
\be T(m\cdot a)=T(m)\cdot a\quad \rho_{M\p}\circ T=(\id_C\o T)\circ \rho_M
\ee
for all $m\in M, a\in A$.

Similarly, $_A{\cal M}(H)^C$ is the category of finite dimensional 
left weak Doi-Hopf modules 
over the left Doi-Hopf datum $(H,A,C)$.

Let us see  what categories $\cat$ are in our earlier examples:

\vbox{
{\noindent\bf Examples:}

{\bf 1} $\cat$ is equivalent to ${\cal M}_{A\equiv H}$, the category of right 
$H$-modules. The equivalence functor $F:\cat \to {\cal M}_{A}$ is the 
forgetful functor.

{\bf 2} $\cat$ is equivalent to $^{C\equiv H}{\cal M}$, the category of 
left $H$-comodules. The equivalence functor ${\hat F}:\cat\to ^C{\cal M}$
is the forgetful functor.

{\bf 3} $\cat$ is equivalent to $^H{\cal M}_H$, the category of weak Hopf
modules \cite{BSz,BNSz} over $H$.

{\bf 4} $\cat$ is equivalent to ${\cal YD}(K^{op}_{cop})^{op}$, the category 
of (some twisted version of) Yetter-Drinfeld modules over $H$. (For its 
definition see the Appendix).
}

\begin{prop} Let $(H,A,C)$ be a finite dimensional right weak Doi-Hopf 
datum and $({\hat H},
{\hat C}, {\hat A})$ its dual. Then  the categories $\cat$ and 
$_{\hat C}{\cal M}({\hat H})^{\hat A} $ are equivalent. 
\end{prop}

{\em Proof:} Let us define the functor $D:\cat \to _{\hat  C}
{\cal M} ({\hat H})^{\hat A}$ 
\bea D(M)\colon = {\hat M} \ {\rm as\  a}\  k-{\rm space}& \qquad 
&{\hat c}\cdot \mu \colon = ({\hat c}\o \mu)\circ \rho_M \nn
&\qquad& {\hat \rho_{\hat M}}(\mu)\colon = a_i \triangleright 
\mu \o \alpha^i\nn
D(T)\colon = T^t &\qquad &
\eea
where $M$ is an object and $T$ an arrow in $\cat$, $^t$ means 
transposition of linear operators, ${\hat c}\in 
{\hat C}$, $\mu\in {\hat M}, (a\triangleright \mu)(m)= \mu(m\cdot a)$ for 
$a\in A, \mu\in {\hat M}, m\in M$, $\{a_i\}$ is a basis for $A$ and 
$\{\alpha^i\}$ is the dual basis for ${\hat A}$. One checks by direct
calculation that $D$ defines an equivalence functor. 
\hfill\qed

\begin{prop} \label{adj}
Let $(H,A,C)$ be a non-degenerate right weak Doi-Hopf datum.
Then the  forgetful functor $F:\cat\to {\cal M}_A$ has a left 
adjoint and ${\hat  F} :\cat \to  {^C{\cal M}}$ has a right adjoint.
\end{prop}

{\em Proof:} Our proof is consructive. Define $G : {\cal M}_A\to \cat$ by
\bea &G(M)\colon =& C\cdot 1_{A<-1>}\o M\cdot 1_{A<0>} \quad  {\rm as\ a}\ k-
{\rm space} \nn
&&(c\o m)\cdot a \colon =
c\cdot a_{<-1>}\o m\cdot a_{<0>} \nn
&&\rho_{G(M)}
\colon = (\Delta_C \o \id_M)\vert_{G(M)} \nn
&G(T)\colon =& (id_C \o T) 
\eea
for $M$ an object and $T$ an arrow in $\cat$, $a\in A, (c\o m)\in G(M)
\subset C\o M$.

The fact that $G$ is a left adjoint of $F$ is justified by the existence of 
unit and counit natural homomorphisms $\rho:\id_{\cat}\to G\circ F$ and
$\delta :F\circ G\to \id_{{\cal M}_A}$. Define them as
\bea
\rho_M : M\to G(M) &\qquad &
\rho_M(m)\colon = m_{<-1>}\o m_{<0>} \\
\delta_M : G(M) \to M &\qquad &
\delta_M \colon = 
(\varepsilon_C\o \id_M)\vert_{G(M)} .
\nonumber\eea
It is staightforward to show that 
$\rho_M\in (M, G(M))_{\cat}$, and $\rho$ 
is natural. The proof of $\delta_M\in (G(M),M)
_{{\cal M}_A}$
lies on the following

\begin{lemma} \label{lemma}
Let $(H,A,C)$ be a non-degenerate right weak Doi-Hopf datum. Then for 
any $c\in C$ and $a\in A$ 
\bea (i)&\ &\Delta_C(c\cdot 1_{A<-1>})\o 1_{A<0>}= 
c_{(1)}\o c_{(2)}\cdot 1_{A<-1>} \o 1_{A<0>} \\
(ii) &\ & \Pi^L(a_{<-1>})\o a_{<0>}= \Pi^L(1_{A<-1>})\o 1_{A<0>}a.
\eea
\end{lemma}

Lemma \ref{lemma} (ii) implies $\varepsilon_C (c\cdot a_{<-1>}) a_{<0>}=
\varepsilon_C(c\cdot 1_{A<-1>})1_{<0>}a$ and hence $\delta_M\in
(G(M),M)_{{\cal M}_A}$. Naturality of 
$\delta$ is obvious.

One can proceed the same way in the case of ${\hat F}$ using now Lemma
\ref{lemma} (i). Define ${\hat G}:
^C{\cal M}\to \cat $ as
\bea {\hat G}(M)&\colon =& 
\{ \varepsilon_C (m_{<-1>}\cdot a_{<-1>}) m_{<0>}\o a_{<0>} 
\vert m\in M, a\in A \}\  {\rm as \ a}\  k-{\rm space} \nn
&&(m\o a)\cdot b \colon =\varepsilon_C(m_{<-1>}\cdot a_{<-1>} b_{<-1>}) 
m_{<0>}\o a_{<0>} b_{<0>} \nn
&& \rho_{{\hat G}(M)}(m\o a) \colon = m_{<-1>}\cdot a_{<-1>}\o m_{<0>} \o
a_{<0>} \nn
{\hat G}(T)&\colon = &T\o \id_A 
\eea
for $M$ an object and $T$ an arrow in $\cat$, $(m\o a)\in {\hat G}(M)
\subset M\o A, b\in A$.

The unit and counit natural homomorphisms ${\hat \rho }:\id_{^C{\cal M}}\to
{\hat F}\circ {\hat G}$ and ${\hat \delta}: {\hat G}\circ {\hat F}\to
\id_{\cat}$ can be given by
\bea {\hat \rho}_M : M\to {\hat G}(M) &\qquad &
{\hat \rho}_M(m)\colon = \varepsilon_C (m_{<-1>}\cdot 1_{A<-1>}) m_{<0>}\o
1_{A<0>} \nn
{\hat \delta}_M: {\hat G}(M)\to M &\qquad &
{\hat \delta}_M(m\o a)\colon = m\cdot a
\eea

\hfill\qed

\sec{The Weak Smash Product}

\begin{defi} For the non-degenerate right weak Doi-Hopf datum $(H,A,C)$ 
define the weak smash 
product algebra $A\#{\hat C}$ as the $k$-space $1_{A<0>}A\o 1_{A<-1>}
\triangleright {\hat C}$  equipped with the multiplication rule 
\be (a\#{\hat c})(b\#{\hat d})\colon = (a_{<0>}b \# 
{\hat c}(a_{<-1>} \triangleright {\hat d})) \label{smash} \ee
for $(a\# {\hat c}), (b\#{\hat d})\in A\#{\hat C}$.
\end{defi}
One checks that (\ref{smash}) makes $A\#{\hat C}$ an associative algebra with 
unit element $1_{A<0>}\#$ $ 1_{A<-1>}\triangleright 1_{\hat C}$.

Let us see what algebras  $A\#{\hat C}$ are  in our earlier examples.

{\noindent \bf Examples:}

{\bf 1} $(A\equiv H) \#({\hat C}\equiv {\hat H}^R)$ is isomorphic to $H$, 
the isomorphism being given by $\iota:A\# {\hat C}\to H$, $\iota\colon =
\id_H\o \varepsilon_{\hat H}$.

{\bf 2} $(A\equiv H^L) \# ({\hat C}\equiv {\hat H})$ is isomorphic to 
${\hat H}$,
the isomorphism being given by $\iota : A\#{\hat C}\to {\hat H}$, $\iota
\colon = \varepsilon \o \id_{\hat H}$.

{\bf 3} $(A\equiv H) \# ({\hat C}\equiv {\hat H})$ is isomorphic to the Weyl 
algebra or Heisenberg double ${\hat H}\cp H$ \cite{BSz,BNSz}, the 
isomorphism being given by $\iota : A\#{\hat C}\to {\hat H}\cp H$,
$\iota (\1_{(2)}a \# (\1_{(1)}\rightharpoonup {\phi}))\colon = \phi a$.
(In all of the examples $\rightharpoonup$ denotes the Sweedler's arrow 
\cite{Sw}.)

{\bf 4} $(A\equiv K) \# ({\hat C}\equiv {\hat K})$ is isomorphic to 
the (twisted)
Drinfel'd double ${\cal D} (K^{op}_{cop})^{op}$ (for its definition see 
the Appendix). The equivalence is given by $\iota: A\#{\hat C}\to {\cal D}
(K^{op}_{cop})^{op}$, $\iota(\1_{(2)}a\# (\1_{(1)}\rightharpoonup {\phi}
\leftharpoonup 
S\inv(\1_{(3)})))\colon = {\cal D}(\phi){\cal D}(a)$.

\begin{prop} \label{cateq}
Let $(H,A,C)$ be a non-degenerate right weak Doi-Hopf datum such that $C$ is 
finite dimensional as a $k$-space. Then the categories $\cat$ and 
${\cal M}_{A\# {\hat C}}$ are isomorphic.
\end{prop}

{\em Proof:} We have the functor $P: \cat \to {\cal M}_{A\#{\hat C}}$
\bea P(M)\colon =M \ {\rm as\ a }\ k-{\rm space} &\quad &
m\cdot (a\#{\hat c})\colon = {\hat c}(m_{<-1>})m_{<0>}\cdot a \nn
P(T)\colon = T &\quad & 
\eea
for $M$ an object and $T$ an arrow in $\cat$, $(a\# {\hat c})\in 
A\#{\hat C}, m\in M$.

If $C$ is finite dimensional as a $k$-space then let $\{c_i\}$ be any basis 
for $C$ and $\{\gamma^i\}$ the dual basis for ${\hat C}$ and construct the 
inverse functor  $P\p:{\cal M}_{A\#{\hat C}}\to \cat$ of $P$:
\bea 
P\p(M)\colon = M \ {\rm as \ a }\ k-{\rm space}&\quad &
m\cdot a\colon = m\cdot (1_{A<0>}a\# 1_{A<-1>}\triangleright 1_{\hat C})\nn
&\quad & \rho_M(m)\colon = c_i \o m\cdot (1_{A<0>}\# 1_{A<-1>}\triangleright 
\gamma^i)\nn
P\p(T)\colon =T &\quad&
\eea
for $M$ an object and $T$ an arrow of $\cat$, $a\in A, m\in M$.

\hfill\qed

\sec{Integrals for Weak Doi-Hopf Data}

Let $(H,A,C)$ be a non-degenerate right weak Doi-Hopf datum where $H$ is a 
weak {\em Hopf} algebra with antipode $S$, $F:\cat \to {\cal M}_A$
the forgetful functor, $G$ its left adjoint as in Proposition \ref{adj}.
$V$ be the $k$-space of the natural homomorphisms $\nu: G\circ F \to 
\id_{\cat}$, called the {\em space of integrals} for the weak Doi-Hopf datum 
$(H,A,C)$. We have a straightforward generalization of Theorem 2.3 
of \cite{CIMZ}:

\begin{theorem} The space $V$ is isomorphic to the space $V_4$:
\bea V_4\colon =\{ \gamma : C\to (C,A)_{\rm Lin} &\vert & 
\forall c,d\in C \quad a\in A \nn
\gamma (c)(d)a &=&
a_{<0>} \gamma (c\cdot a_{<-2>})(d\cdot a_{<-1>}) \nn
c_{(1)} \o \gamma(c_{(2)})(d)&=& d_{(2)}\cdot \gamma(c)(d_{(1)})_{<-1>}\o
\gamma(c)(d_{(1)})_{<0>} \}. 
\eea
Furthermore the isomorphism $f_4:V\to V_4$
takes $\nu\in V$ to a {\em normalized element} of $V_4$ i.e.
to an element $\gamma\in V_4$ such that $\gamma(c_{(1)})(c_{(2)})=
\varepsilon_C (c\cdot 1_{A<-1>}) 1_{A<0>}$ 
if and only if $\nu$ is a splitting of the unit natural homomorphism 
$\rho:\id_{\cat}\to G\circ F$. 
\end{theorem}
The relevance of the existence of normalized elements in $V_4$ is discussed 
in \cite{CIMZ}.

Let us turn to the investigation of the space of integrals over the 
weak Doi-Hopf datum $(H,A,C)$ in our earlier examples. In doing so
we make the additional assumption in the Examples {\bf 2} and {\bf 3} 
on $H$ and in {\bf 4} on $K$ to be a 
{\em Frobenius} WHA. Under this additional condition we identify the space 
of integrals for the weak Doi-Hopf datum $(H,A,C)$ with certain subspace
of the smash product algebra $A\#{\hat C}$. Also the normalization condition 
is formulated as a relation in the algebra $A\#{\hat C}$.

In all of the examples $r$ be a non-degenerate right integral in $H$
and $\rho$ the dual right integral \cite{BNSz} in ${\hat H}$.

{\noindent \bf Examples:}

{\bf 1} The space of Doi-Hopf integrals over $(H,A,C)$ 
is isomorphic to $V_0\colon = {\rm Center} H$.
Construct the isomorphism  $f:V_4\to V_0$ as 
\be f(\gamma)\colon =\gamma(\1)(\1).\ee
The unique normalized element of $V_0$ is the unit element $\1$ of $H$.

{\bf 2} The space of the Doi-Hopf integrals is isomorphic to
$V_0\colon = ({\hat H}^{R})\p\cap {\hat H}$, the commutant of the right
subalgebra in ${\hat H}$. Let us construct the isomorphism $f:V_4\to V_0$ as 
\be [f(\gamma)](h)\colon = \varepsilon(\gamma(r)(h))\ee
for all $h\in H$. 

An element $\xi\in V_0$ is normalized if 
\be {\hat S}\inv(\rho_{(2)})\xi \rho_{(1)} ={\hat \1} \ee
holds in ${\hat H}$.

The space $V_0$ is {\em not} isomorphic to 
the space ${\cal I}^L({\hat H})$ of left integrals in ${\hat H}$. It is its 
subspace ${\hat H}^L$ which is isomorphic to ${\cal I}^L({\hat H})$ via the 
isomorphism $g: {\cal I}^L({\hat H})\to {\hat H}^L$,
$g(\lambda)\colon = {\hat S} (\lambda\leftharpoonup r)$. It is but true 
that the existence of normalized elements in ${\cal I}^L({\hat H})$ and 
$V_0$ are equivalent. 

{\bf 3} The space of the Doi-Hopf integrals is isomorphic to
$V_0\colon = H\p\cap ({\hat H}\cp H)$, 
the commutant of $H$ in the Weyl algebra.
The  isomorphism $f: V_4\to V_0$ is given by 
\be f(\gamma)\colon = \beta^i \gamma(r)(b_i) \ee
with the help of the basis $\{b_i\}$ of $H$ and the dual basis $\{\beta^i\}$
of ${\hat H}$.

The element $w\in V_0$ is normalized if
\be {\hat S}\inv (\rho_{(2)})w\rho_{(1)}=1_{{\hat H}\cp H}  \ee
holds in the Weyl algebra ${\hat H}\cp H$.

{\bf 4} The space of the Doi-Hopf integrals is isomorphic to
$V_0\colon =\{ u\in {\cal D}(K^{op}_{cop})^{op}
\vert$ $ u {\cal D}(b)\ =\ {\cal D}(b_{(1)}) u 
{\cal D}( S\inv (r) S^{-2}(b_{(2)})
\rightharpoonup \rho) \}$. The isomorphism 
$f:V_4\to V_0$ is given by
\be f(\gamma)\colon = {\cal D}(\beta^i) {\cal D}(\gamma(r)(b_i)) \ee
with the help of the basis $\{b_i\}$ of $H$ and the dual basis $\{\beta^i\}$
of ${\hat H}$.

$u\in V_0 $ is normalized if 
\be {\hat S}\inv (\rho_{(2)})u\rho_{(1)}=1_{{\cal D}(K^{op}_{cop})^{op}}  \ee
holds in the double ${\cal D}(K^{op}_{cop})^{op}$.

\sec{Appendix: Yetter-Drinfel'd modules over WHA's and Drinfel'd doubles}

For the convenience of the reader we give here 
the generalization of the double construction due to Drinfel'd 
\cite{Dr} and of the corresponding theory of Yetter-Drinfel'd modules 
\cite{Y,RT} to WHA's.

\begin{defi} \cite{BSz} Let $H$ be a finite dimensional WHA over the field $k$.
Its {\em Drinfel'd double} ${\cal D}(H)$ is the WHA defined below:

As a $k$-space ${\cal D}(H)$ is an amalgamated tensor product 
$H_{H^L\equiv {\hat H}^R}\o_{H^R\equiv {\hat H}^L} {\hat H}$ with 
the amalgamation relations $a^R\o {\hat \1}\equiv \1\o ({\hat \1}\leftharpoonup 
a^R)$; $(a^L\o {\hat \1})\equiv \1\o (a^L\rightharpoonup {\hat \1})$ for $a^L\in H^L, 
a^R\in H^R$. Denote by ${\cal D}(a){\cal D}(\phi)$ the image of 
$H\o {\hat H}\ni a\o\phi$ under the amalgamation and ${\cal D}(a)\equiv
{\cal D}(a){\cal D}({\hat \1}), {\cal D}(\phi)\equiv {\cal D}(\1){\cal D}
(\phi)$.

The algebra structure is defined by 
\bea {\cal D}(a) {\cal D}(b) &=& {\cal D} (ab) \nn
     {\cal D}(\phi){\cal D}(\psi) &=& {\cal D}(\phi\psi) \nn
     {\cal D}(\phi) {\cal D}(a) &=& {\cal D} (a_{(2)}) {\cal D}(\phi_{(2)})
\pair{\phi_{(1)}}{a_{(3)}}\pair{\phi_{(3)}}{S\inv(a_{(1)})}.
\label{dbalg}\eea

One checks that (\ref{dbalg}) is compatible with the amalgamation relations 
and makes ${\cal D}(H)$ an associative algebra with unit ${\cal D}(\1)
\equiv {\cal D}({\hat \1})$.

The colagebra structure is given by
\bea \Delta_{\cal D}({\cal D}(a){\cal D}(\phi)) &=& 
{\cal D}(a_{(1)}) {\cal D} (\phi_{(2)}) \o
{\cal D}(a_{(2)}) {\cal D} (\phi_{(1)}) \nn
\varepsilon_{\cal D} ({\cal D}(a) {\cal D}(\phi)) &=&
\varepsilon (a (\phi \rightharpoonup \1))\equiv 
{\hat \varepsilon}(({\hat \1}\leftharpoonup a)\phi ).
\label{dbcoalg}\eea
One checks that (\ref{dbcoalg}) makes ${\cal D}(H)$ a WBA. Finally the 
antipode is
\be S_{\cal D}({\cal D}(a) {\cal D}(\phi))={\cal D}({\hat S}\inv (\phi)) 
{\cal D}(S(a)) \ee
making ${\cal D}(H)$ a WHA.
\end{defi}

\begin{defi} Let $H$ be  a WBA over the field $k$. The $k$-space $M$ 
is a right {\em Yetter-Drinfel'd module} over $H$ if it is a non-degenerate 
right $H$-module and a non-degenerate left $H$ comodule s.t. 
\bea m_{<-1>} a_{(1)} \o m_{<0>}\cdot a_{(2)} &=&
a_{(2)} (m\cdot a_{(1)})_{<-1>}\o (m\cdot a_{(1)})_{<0>} \nn
m_{<-1>} \1_{(1)}\o m_{<0>}\cdot \1_{(2)} &=& m_{<-1>} \o m_{<0>}
\label{YD}\eea
for all $m\in M,a\in A$.
\end{defi}

Notice that if $H$ is also a WHA then (\ref{YD}) can be replaced by the
single relation
\be (m\cdot a)_{<-1>}\o (m\cdot a)_{<0>}=
S\inv(a_{(3)})m_{<-1>}a_{(1)} \o m_{<0>}\cdot a_{(2)}. \ee

By the category ${\cal YD}(H)$ we mean the category with objects the
finite dimensional right Yetter-Drinfel'd modules over $H$ and 
arrows $T:M\to M\p$ intertwining
both the $H$-module and the $H$-comodule structures of $M$ and $M\p$.

If $H$ is a finite dimensional WHA then by our Proposition  \ref{cateq} 
and Example 4. the
category ${\cal YD}(H)$ is equivalent to the category of the right modules
over the WHA  ${\cal D}(H)$ hence 
carries (among others) a monoidal structure \cite{BSz,BNSz}. It is not 
so obvious however that it is true for any WBA $H$:

\begin{prop}  Let $H$ be a WBA over the field $k$. Then the category 
${\cal YD}(H)$ has a monoidal structure.
\end{prop}

{\em Proof:} Our proof is constructive. For two objects $M,N$ and arrows
$T,S$ of  ${\cal YD}(H)$ let
\bea &M\times N &\colon = M\cdot \1_{(1)} \o N\cdot \1_{(2)}
\quad {\rm as\ a}\ k-{\rm space} \nn
&&(m\o n)\cdot a\colon = m\cdot a_{(1)}\o n\cdot a_{(2)} \nn
&&\rho_{M\times N}(m\o n)\colon = n_{<-1>} m_{<-1>}\1_{(1)}\o 
m_{<0>}\cdot \1_{(2)}\o n_{<0>}\cdot \1_{(3)} \nn
&T\times S&\colon = (T\o S)\circ \Delta(\1) 
\eea
with $m\o n\in M\times N, a\in A$. The monoidal unit is 
\bea H^L\ {\rm as\ a}\ k-{\rm space} &\quad& a^L\cdot h \colon = 
\1_{(2)} \varepsilon(a^Lh\1_{(1)}) \nn
&\quad& \rho_{H^L} \colon =\Delta\vert_{H^L}
\eea
for $a^L\in H^L, h\in H$. 

The reader may check using some WBA calculus that all $M\times N$ and $H^L$
are Yetter-Drinfel'd modules over $H$ if $M$ and $N$ are. 

In order to prove that $H^L$ is a monoidal unit for the category 
${\cal YD}(H)$ one has to construct the invertible intertwiners $u_M^L\in
(M,H^L\times M)_{{\cal YD}(H)}$, $u_M^R\in (M, M\times H^L)_{{\cal YD}(H)}$ 
satisfying the triangle identities \cite{McL} and being natural in $M$. 
They are as follows:
\bea u_M^L(m) &=& \1_{(2)}\o m\cdot \Pi^L(\1_{(1)}) \nn
     u_M^R(m) &=& m\cdot \1_{(1)} \o \1_{(2)}.
\eea
for all $m\in M$ and all objects $M$ of ${\cal YD}(H)$.
\hfill\qed


\begin{thebibliography}{FroGab}
\addcontentsline{toc}{section}{\protect\numberline{}{References}}
\renewcommand{\b}{\bibitem}
\renewcommand{\baselinestretch}{.3}
\small
\medskip


\b{BSz} G. B\"ohm, K. Szlach\'anyi, {\em Lett. Math. Phys} {\bf 35} (1996)
p.437 
\b{BNSz} G. B\"ohm, F. Nill, K. Szlach\'anyi, {\em `Weak Hopf Algebras I: 
Integral Theory and $C^*$-structure'}  math.QA/9805116  to appear in {\em
J. Algebra}
\b{BNSz2} G. B\"ohm, F. Nill, K. Szlach\'anyi, {\em `Weak Hopf Algebras II:
Representation Theory, Dimensions, and the Markov Trace'} in preparation
\b{CMZ1} C. Caenepeel, G. Militaru, S. Zhu, {\em Israel J. Math} {\bf 100} 
(1997) p.221
\b{CMZ2}  C. Caenepeel, G. Militaru, S. Zhu, {\em Trans AMS} {\bf 349} (1997)
p.4311
\b{CIMZ} C. Caenepeel, G. Militaru, S. Zhu, {\em J. Algebra} {\bf 187} (1997) 
p.388;  C. Caenepeel, B. Ion, G. Militaru, S. Zhu {\em `Separable functors
for the category of Doi-Hopf modules. Applications'} math.QA/9809021
\b{Doi} Y. Doi, {\em J. Algebra} {\bf 153} (1992) p.373
\b{Doint} Y. Doi, {\em Comm. Algebra} {\bf 13} (1985) p.2137
\b{Dr} V. G. Drinfel'd, {\em Proc. Int. Congr. Math. Berkeley} (1986) p.798 
\b{Kop} M. Koppinen, {\em `Variations of the smash product with applications
to group graded rings'} preprint 1991
\b{McL} S. McLane, {\em `Categories for the working mathematician'} 
Springer 1971
\b{N} F. Nill, {\em `Axioms for Weak Bialgebras'} math.QA/9805104
\b{NSzW} F. Nill, K. Szlach\'anyi, H-W Wiesbrock, {\em `Weak Hopf Algebras 
and Reducible Jones Inclusions of Depth 2'} math.QA/9806130
\b{RT} D. Radford, J. Towber, {\em J. Pure and Appl. Algebra} {\bf 87} (1993)
p.259
\b{Sw} M. E. Sweedler, {\em `Hopf algebras'} Benjamin 1969
\b{Y} D. N. Yetter, {\em Math. Proc. Cambridge Phil. Soc.} {\bf 108} (1990)
p.261

\end{thebibliography}
\end{document}